\documentclass[12pt]{article}
\usepackage{amssymb, amsmath}
\usepackage{float,epsfig}

\begin{document}

\newtheorem{theorem}{\bf Theorem}[section]
\newtheorem{proposition}[theorem]{\bf Proposition}
\newtheorem{definition}[theorem]{\bf Definition}
\newtheorem{corollary}[theorem]{\bf Corollary}
\newtheorem{example}[theorem]{\bf Example}
\newtheorem{exam}[theorem]{\bf Example}
\newtheorem{remark}[theorem]{\bf Remark}
\newtheorem{lemma}[theorem]{\bf Lemma}
\newtheorem{conclusions}[theorem]{\bf Conclusions}
\newcommand{\nrm}[1]{|\!|\!| {#1} |\!|\!|}

\newcommand{\ba}{\begin{array}}
\newcommand{\ea}{\end{array}}
\newcommand{\von}{\vskip 1ex}
\newcommand{\vone}{\vskip 2ex}
\newcommand{\vtwo}{\vskip 4ex}
\newcommand{\dm}[1]{ {\displaystyle{#1} } }

\newcommand{\be}{\begin{equation}}
\newcommand{\ee}{\end{equation}}
\newcommand{\beano}{\begin{eqnarray*}}
\newcommand{\eeano}{\end{eqnarray*}}
\newcommand{\inp}[2]{\langle {#1} ,\,{#2} \rangle}
\def\bmatrix#1{\left[ \begin{matrix} #1 \end{matrix} \right]}
\def \noin{\noindent}
\newcommand{\evenindex}{\Pi_e}



\def \R{{\mathbb R}}
\def \C{{\mathbb C}}
\def \K{{\mathbb K}}
\def \J{{\mathbb J}}
\def \Lb{\mathrm{L}}

\def \T{{\mathbb T}}
\def \Pb{\mathrm{P}}
\def \N{{\mathbb N}}
\def \Ib{\mathrm{I}}
\def \Ls{{\Lambda}_{m-1}}
\def \Gb{\mathrm{G}}
\def \Hb{\mathrm{H}}
\def \Lam{{\Lambda_{m}}}
\def \Qb{\mathrm{Q}}
\def \Rb{\mathrm{R}}
\def \Mb{\mathrm{M}}
\def \norm{\nrm{\cdot}\equiv \nrm{\cdot}}

\def \P{{\mathbb P}_m(\C^{n\times n})}
\def \A{{{\mathbb P}_1(\C^{n\times n})}}
\def \H{{\mathbb H}}
\def \L{{\mathbb L}}
\def \G{{\mathcal G}}
\def \S{{\mathbb S}}
\def \sigmin{\sigma_{\min}}
\def \elam{\sigma_{\epsilon}}
\def \slam{\sigma^{\S}_{\epsilon}}
\def \Ib{\mathrm{I}}
\def \Tb{\mathrm{T}}
\def \d{{\delta}}

\def \Lb{\mathrm{L}}
\def \N{{\mathbb N}}
\def \Ls{{\Lambda}_{m-1}}
\def \Gb{\mathrm{G}}
\def \Hb{\mathrm{H}}
\def \Delta{\triangle}
\def \Rar{\Rightarrow}
\def \p{{\mathsf{p}(\lam; v)}}

\def \D{{\mathbb D}}

\def \tr{\mathrm{Tr}}
\def \cond{\mathrm{cond}}
\def \lam{\lambda}
\def \sig{\sigma}
\def \sign{\mathrm{sign}}

\def \ep{\epsilon}
\def \diag{\mathrm{diag}}
\def \rev{\mathrm{rev}}
\def \vec{\mathrm{vec}}

\def \spsd{\mathsf{spsd}}
\def \spd{\mathsf{spd}}
\def \sk{\mathsf{skew}}
\def \sy{\mathsf{sym}}
\def \en{\mathrm{even}}
\def \odd{\mathrm{odd}}
\def \rank{\mathrm{rank}}
\def \pf{{\bf Proof: }}
\def \dist{\mathrm{dist}}
\def \rar{\rightarrow}

\def \rank{\mathrm{rank}}
\def \pf{{\bf Proof: }}
\def \dist{\mathrm{dist}}
\def \Re{\mathsf{Re}}
\def \Im{\mathsf{Im}}
\def \re{\mathsf{re}}
\def \im{\mathsf{im}}

\def \sym{\mathsf{sym}}
\def \sksym{\mathsf{skew\mbox{-}sym}}
\def \odd{\mathrm{odd}}
\def \even{\mathrm{even}}
\def \herm{\mathsf{Herm}}
\def \skherm{\mathsf{skew\mbox{-}Herm}}
\def \str{\mathrm{ Struct}}
\def \eproof{$\blacksquare$}
\def \proof{\noin\pf}

\def \bS{{\bf S}}
\def \cA{{\cal A}}
\def \E{{\mathcal E}}
\def \X{{\mathcal X}}
\def \F{{\mathcal F}}
\def \tr{\mathrm{Tr}}
\def \range{\mathrm{Range}}

\def \pal{\mathrm{palindromic}}
\def \palpen{\mathrm{palindromic~~ pencil}}
\def \palpoly{\mathrm{palindromic~~ polynomial}}
\def \odd{\mathrm{odd}}
\def \even{\mathrm{even}}


\title{Vector space of linearizations for the quadratic two-parameter matrix polynomial}

\author{ Bibhas Adhikari\thanks{E-mail:
bibhas@iitj.ac.in} \\ Indian Institute of Technology Rajasthan \\ Jodhpur, India }
\date{}


\maketitle \thispagestyle{empty}

\noin\textbf{Abstract:} Given a quadratic two-parameter matrix polynomial $Q(\lam,\mu)$, we develop a systematic approach to generating a vector space of linear two-parameter matrix polynomials. We identify a set of linearizations of $Q(\lam,\mu)$ that lie in the vector space. Finally, we determine a class of linearizations for a  quadratic two-parameter eigenvalue problem.\\

\noin\textbf{Key words:} matrix polynomial, two-parameter matrix polynomial, quadratic two-parameter eigenvalue problem, two-parameter
eigenvalue problem, linearization \\

\noin\textbf{AMS classification:} 65F15, 15A18, 15A69, 15A22

 \section{Introduction}

We consider two-parameter quadratic matrix polynomials of the form \begin{equation}\label{def:qtp} Q(\lam, \mu) =A + \lam B + \mu C + \lam^2 D + \lam\mu E + \mu^2 F \end{equation} where $\lam, \mu$ are scalars and the coefficient matrices are real or complex matrices of order $n\times n.$
If $(\lam,\mu)\in \C\times\C$ and nonzero $x\in \C^n$ satisfy $Q(\lam,\mu)x = 0$, then $x$ is said
to be an eigenvector of $Q(\lam, \mu)$ corresponding to the eigenvalue $(\lam,\mu)$.
The classical approach to solving spectral problems for matrix polynomials
is to first perform a linearization, that is, to transform the given
polynomial into a linear matrix polynomial, and then work with this linear polynomial (see \cite{Kha97,Kub98,Kha07,mupl08,mupl081,mackey3} and the references therein). Therefore, given a quadratic two-parameter matrix polynomial $Q(\lam,\mu)$, we seek linear two-parameter matrix polynomials $L(\lam,\mu)= \lam \widehat{A}_1 + \mu \widehat{A}_2 + \widehat{A}_3$, called \emph{linearizations}, which have the same eigenvalues as $Q(\lam,\mu).$

The one-parameter matrix polynomials have been an active area of research in numerical linear algebra \cite{mackey3,fiedler1,fiedler2}. In \cite{mackey3}, Mackey et al. have investigated the one-parameter polynomial eigenvalue problem extensively and they have produced vector spaces of linearizations for a one-parameter matrix polynomial by generalizing the companion forms of the one-parameter polynomial. Adopting a similar approach we derive a set of linearizations of a quadratic two-parameter matrix polynomial.

%

 The quadratic two-parameter eigenvalue problem is concerned with finding scalars $\lam, \mu \in\C$ and non-zero vectors $x_1\in\C^{n_1}, x_2\in\C^{n_2}$ such that \begin{equation}\label{defn:qevp1} \left\{
\begin{array}{ll}
 Q_1(\lam,\mu)x_1=(A_1 + \lam B_1 + \mu C_1 + \lam^2 D_1 + \lam\mu E_1 + \mu^2 F_1)x_1=0 & \hbox{} \\
 Q_2(\lam,\mu)x_2=(A_2 + \lam B_2 + \mu C_2 + \lam^2 D_2 + \lam\mu E_2 + \mu^2 F_2)x_2=0 & \hbox{}
  \end{array}
\right.
  \end{equation} where $A_i, B_i, \hdots, F_i\in\C^{n_i\times n_i}, i=1,2.$ A pair $(\lam, \mu)$ satisfying (\ref{defn:qevp1}) is called an eigenvalue of (\ref{defn:qevp1}) and $x_1\otimes x_2,$ where $\otimes$ is the Kronecker product, is the corresponding eigenvector. This problem appears in stability analysis of different systems, for example, time-delay systems of single delay \cite{HocMuPl10,jar08,jarhocs09,mupl08}. The standard approach to solving (\ref{defn:qevp1}) is by linearizing the problem into a two-parameter eigenvalue problem of larger size and then by converting it into an equivalent coupled generalized eigenvalue problem which is then solved by numerical methods, see \cite{mupl08,mupl081,HocMuPl10}.

  Given (\ref{defn:qevp1}), we seek a two-parameter  eigenvalue problem \begin{equation}\label{sec1:defn:qlin} \left\{
                     \begin{array}{ll}
                       L_1(\lam,\mu)w_1 := (A^{(1)}) + \lam B^{(1)} + \mu C^{(1)})w_1 =0& \hbox{} \\
                       L_2(\lam,\mu)w_2 := (A^{(2)}) + \lam B^{(2)} + \mu C^{(2)})w_2 =0& \hbox{}
                     \end{array}
                   \right.
   \end{equation} with the same eigenvalues, where $w_i \in\C^{3n_i}\setminus\{0\}$ and $A^{(i)}, B^{(i)}, C^{(i)}\in \C^{3n_i \times 3n_i},$ $ i=1,2.$ In such case (\ref{sec1:defn:qlin})   is called a linearization of (\ref{defn:qevp1}).


  The choice of a linearization may have an adverse effect on the sensitivity of the eigenvalues. Therefore, it is important to identify potential linearizations and describe their constructions. In this paper, we develop a systematic approach that enables us to generate a class of linearizations for a quadratic two-parameter eigenvalue problem.

\section{Linearizations for quadratic two-parameter matrix polynomial}

 In this section we construct a set of linearizations of a quadratic two-parameter matrix polynomial.

 \begin{definition}(\cite{mupl08})
 A $ln\times ln$ linear matrix polynomial $L(\lam, \mu) = \lam \widehat{A}_1 + \mu \widehat{A}_2 + \widehat{A}_3$ is a linearization of an $n\times n$ matrix polynomial $Q(\lam, \mu)$ if there exist polynomials $P(\lam, \mu)$ and $R(\lam, \mu),$ whose determinant is a non-zero constant independent of $\lam$ and $\mu,$ such that $$\bmatrix{Q(\lam,\mu) & 0 \\ 0 & I_{(l-1)n}} = P(\lam, \mu)L(\lam, \mu)R(\lam, \mu).$$

 \end{definition}

 Let $Q(\lam, \mu)$ be a quadratic two-parameter matrix polynomial given by $$Q(\lam, \mu)=\lam^2A_{20} + \mu^2A_{02}+\lam\mu A_{11}+\lam A_{10}+\mu A_{01}+A_{00}$$ where the coefficient matrices are of order $n\times n.$ Assume that $x$ is the eigenvector corresponding to an eigenvalue $(\lam,\mu)$ of $Q(\lam, \mu),$ that is, $Q(\lam, \mu)x =0.$ With a view to constructing linearizations of $Q(\lam, \mu)$, we denote $x=x_{00}, \lam x_{00}= x_{10}, \mu x_{00} = x_{01}.$ Then we have \begin{equation} A_{20}(\lam x_{10}) + A_{02}(\mu x_{01}) + A_{11} (\lam x_{01}) + A_{10}x_{10} + A_{01}x_{01} + A_{00}x_{00} =0. \end{equation}


%
%
  Consequently we have \begin{equation}\label{def:l1} \left(\lam \underbrace{\bmatrix{A_{20} & A_{11} & 0 \\ 0 & 0 & 0 \\ 0 & 0 & I}}_{\widehat{A}_1} + \mu \underbrace{\bmatrix{0 & A_{02} & 0 \\ 0 & 0 & I\\ 0 & 0 & 0}}_{\widehat{A}_2} + \underbrace{\bmatrix{A_{10} & A_{01} & A_{00} \\ 0 & -I & 0\\ -I & 0 & 0}}_{\widehat{A}_3} \right) \bmatrix{x_{10} \\ x_{01} \\ x_{00}} = 0.\end{equation} Observe that $\bmatrix{x_{10} \\ x_{01} \\ x_{00}} = \bmatrix{\lam x \\ \mu x \\ x} = \bmatrix{\lam\\ \mu \\ 1}\otimes x.$ We denote $\Lambda:= \bmatrix{\lam\\ \mu \\ 1}.$ Thus $x$ is the eigenvector corresponding to an eigenvalue $(\lam, \mu)$ of $Q(\lam, \mu)$ if and only if $L(\lam, \mu) w =0$ where $w =\Lambda\otimes x$ and $L(\lam, \mu) = \lam \widehat{A}_1 + \mu \widehat{A}_2 + \widehat{A}_3,$ that is, $w$ is the eigenvector corresponding to an eigenvalue $(\lam,\mu)$ of $L(\lam,\mu).$ We show that $L(\lam,\mu)$ is a linearization of $Q(\lam,\mu).$

  Define \beano E(\lam,\mu) &:=& \bmatrix{\lam I_n & I_n & 0\\ \mu I_n & 0 & I_n \\ I_n & 0 & 0},\\ F(\lam,\mu) &:=& \bmatrix{I_n & \mu A_{02}+\lam A_{11}+A_{01} & \lam A_{20}+A_{10} \\ 0 & 0 & I_n \\ 0 & I_n & 0}.\eeano Notice that $E, F$ are unimodular quadratic two-parameter matrix polynomials, that is, determinants of $E$ and $F$ are constants. Then we can easily check that \begin{equation} F(\lam,\mu)L(\lam,\mu)E(\lam,\mu) = \bmatrix{Q(\lam,\mu) & 0\\ 0 & I_{2n}}.\end{equation} Thus we have $\mbox{det}Q(\lam,\mu)=\gamma\mbox{det}L(\lam,\mu)$ for some $\gamma\neq 0.$ This implies $L(\lam,\mu)$ preserves the eigenvalues of $Q(\lam, \mu)$ and hence is a linearization of order $3n\times 3n.$ We call this linearization the \textit{standard} linearization of $Q(\lam,\mu)$. It is interesting to observe that the linearization proposed in \cite{mupl08} is up to some permutations of block rows and columns of the standard linearization.


 However, for the \textit{standard} linearization we have \begin{equation}\label{erual:LQ} L(\lam, \mu) \cdot (\Lambda\otimes x) =\left[ (Q(\lam, \mu)x)^T \,\, 0 \,\, \hdots \,\, 0\right]^T \,\, \mbox{for all} \,\, x\in\C^n,\end{equation} and therefore, any solution of (\ref{def:l1}) gives a solution of the original problem $Q(\lam, \mu)x=0.$ Further, by (\ref{erual:LQ}) we have \begin{equation}\label{lincondition} L(\lam, \mu) \cdot (\Lambda\otimes I_n) = L_1(\lam, \mu) \bmatrix{\lam I_n \\ \mu I_n \\  I_n} = \bmatrix{Q(\lam, \mu) \\ 0 \\ 0} = e_1 \otimes Q(\lam, \mu), e_1=\bmatrix{1\\ 0\\ 0}.  \end{equation}

We restrict our attention to the equation (\ref{lincondition}) which is satisfied by the \textit{standard} linearization. It would be worthy to find linear two-parameter
matrix polynomials $L(\lam,\mu)$ that satisfy  \begin{equation}\label{intro:v}L(\lam,\mu)\cdot (\Lambda\otimes I_n) = v \otimes Q(\lam, \mu)\end{equation} for some vector $v\in\C^3.$ Therefore, we introduce the notation \begin{equation} \mathcal{V}_Q =\{ v\otimes Q(\lam, \mu) : v\in\K^3\} \end{equation} and define \begin{equation}\label{def:L1} \L(Q(\lam, \mu)) := \left\{ L(\lam, \mu) = \lam \widehat{A}_1 + \mu \widehat{A}_2 + \widehat{A}_3, \widehat{A}_i\in\K^{3n\times 3n} : L(\lam, \mu) \cdot (\Lambda \otimes I_n) \in\mathcal{V}_Q\right\}.\end{equation} Note that $\L(Q(\lam, \mu))\neq\emptyset$ as the \textit{standard} linearization $L(\lam,\mu)\in \L(Q(\lam, \mu)).$ It is easy to check that $\L(Q(\lam, \mu)) $ is a vector space. If $L(\lam,\mu)\in \L(Q(\lam, \mu)) $ for some $v\in\C^3$ then call $v$ is an \textit{ansatz} vector associated with $L(\lam,\mu)$.  To investigate the structure of each $L(\lam,\mu)\in \L(Q(\lam, \mu)),$ we define a ``box-addition" for three $3n\times 3n$ block matrices as follows.

%


\begin{definition} Let $\widehat{X}, \widehat{Y}, \widehat{Z} \in \C^{3n\times 3n}$ be three block matrices of the form $$\widehat{X}=\bmatrix{X_{11} & X_{12} & X_{13} \\ X_{21} & X_{22} & X_{23} \\ X_{31} & X_{32} & X_{33}}, \widehat{Y} = \bmatrix{Y_{11} & Y_{12} & Y_{13} \\ Y_{21} & Y_{22} & Y_{23} \\ Y_{31} & Y_{32} & Y_{33}}, \widehat{Z}= \bmatrix{Z_{11} & Z_{12} & Z_{13} \\ Z_{21} & Z_{22} & Z_{23} \\ Z_{31} & Z_{32} & Z_{33}}.$$ Define
\beano  \widehat{X} \boxplus \widehat{Y} \boxplus \widehat{Z}  &=& \bmatrix{X_{11} & X_{12} & X_{13} \\ X_{21} & X_{22} & X_{23} \\ X_{31} & X_{32} & X_{33}} \boxplus \bmatrix{Y_{11} & Y_{12} & Y_{13} \\ Y_{21} & Y_{22} & Y_{23} \\ Y_{31} & Y_{32} & Y_{33}} \boxplus \bmatrix{Z_{11} & Z_{12} & Z_{13} \\ Z_{21} & Z_{22} & Z_{23} \\ Z_{31} & Z_{32} & Z_{33}} \\ &=& \bmatrix{X_{11} & X_{12} & 0 & X_{13} & 0 & 0 \\ X_{21} & X_{22} & 0 & X_{23} & 0 & 0 \\ X_{31} & X_{32} & 0 & X_{33} & 0 & 0} + \bmatrix{0 & Y_{11} & Y_{12} & 0 & Y_{13} & 0 \\ 0 & Y_{21} & Y_{22} & 0 & Y_{23} & 0 \\ 0 & Y_{31} & Y_{32} & 0 & Y_{33} & 0} \\ && + \bmatrix{0 & 0 & 0 & Z_{11} & Z_{12} & Z_{13} \\ 0 & 0 & 0 & Z_{21} & Z_{22} & Z_{23} \\ 0 & 0 & 0 & Z_{31} & Z_{32} & Z_{33}} \eeano where $`+$' is the usual matrix addition.
\end{definition}

For the \textit{standard} linearization $L(\lam, \mu)=\lam \widehat{A}_1 + \mu \widehat{A}_2 + \widehat{A}_3\in\L(Q(\lam, \mu)) $ we have
\beano \widehat{A}_1 \boxplus \widehat{A}_2 \boxplus \widehat{A}_3 &=& \bmatrix{A_{20} & A_{11} & 0 \\ 0 & 0 & 0 \\ 0 & 0 & I} \boxplus \bmatrix{0 & A_{02} & 0 \\ 0 & 0 & I\\ 0 & 0 & 0} \boxplus  \bmatrix{A_{10} & A_{01} & A_{00} \\ 0 & -I & 0\\ -I & 0 & 0} \\ &=& \bmatrix{A_{20} & A_{11} & 0 & 0 & 0 & 0\\ 0 & 0 & 0 &0 & 0 & 0\\ 0& 0 & 0 & I & 0 & 0} + \bmatrix{0 & 0 & A_{02} & 0 & 0 & 0\\ 0 & 0 & 0 &0 & I & 0 \\ 0 & 0 & 0 & 0 & 0 & 0 } \\  && \hfill{ + \bmatrix{ 0 & 0 &0 &A_{10} & A_{01} & A_{00} \\ 0 & 0 & 0& 0 & -I & 0\\ 0 & 0 & 0 & -I & 0 & 0}} \\ &=& \bmatrix{A_{20} & A_{11} & A_{02} & A_{10} & A_{01} & A_{00}\\ 0 & 0 & 0 &0 & 0 & 0\\ 0& 0 & 0 & 0 & 0 & 0} \\ &=& e_1 \otimes \bmatrix{A_{20} & A_{11} & A_{02} & A_{10} & A_{01} & A_{00}} .\eeano Thus we have the following lemma.

\begin{lemma}\label{lem:sec2}
Let $Q(\lam, \mu) = \lam^2A_{20} + \mu^2A_{02}+\lam\mu A_{11}+\lam A_{10}+\mu A_{01}+A_{00}$ be a quadratic two-parameter matrix polynomial with real or complex coefficient matrices of order $n\times n$, and $L(\lam, \mu)= \lam \widehat{A}_1 + \mu \widehat{A}_2 + \widehat{A}_3$ a $3n\times 3n$ two-parameter linear matrix polynomial. Then \beano L(\lam, \mu) \cdot (\Lambda\otimes I_n) &=& v\otimes Q(\lam,\mu) \Leftrightarrow \\ && \widehat{A}_1 \boxplus  \widehat{A}_2 \boxplus \widehat{A}_3 = v\otimes \bmatrix{A_{20} & A_{11} & A_{02} & A_{10} & A_{01} & A_{00}}. \eeano

\end{lemma}

\noin\pf Computational and easy to check.


\begin{example}\label{exp1}

Consider a quadratic two-parameter polynomial $$ Q(\lam, \mu) = \lam^2A_{20} + \mu^2A_{02}+\lam\mu A_{11}+\lam A_{10}+\mu A_{01}+A_{00} $$ where $A_{ij}\in \C^{n\times n}$ and \beano L(\lam,\mu) = && \lam\bmatrix{A_{20} & A_{11}+A_{20} & A_{10}+A_{01} \\ A_{20} & A_{00} & 0 \\ 2A_{20} & A_{02} + 2A_{11} & I } \\ && +\mu \bmatrix{-A_{20} & A_{02} & A_{01} \\ A_{11}-A_{00} & A_{02} & 0\\ -A_{02} & 2A_{02} & A_{01}}  + \bmatrix{-A_{01} & 0 & A_{00} \\ A_{01} & A_{01} & A_{00}\\ -I+2A_{10} & A_{01} & 2A_{00}}.\eeano Then $L(\lam,\mu)\in\L(Q(\lam, \mu))$ since $$\widehat{A}_1 \boxplus \widehat{A}_2 \boxplus \widehat{A}_3 = [1 \,\, 1 \,\, 2]^T \otimes \bmatrix{A_{20} & A_{11} & A_{02} & A_{10} & A_{01} & A_{00}}.$$

\end{example}

Using Lemma \ref{lem:sec2} we characterize the structure of any $L(\lam, \mu) \in \L(Q(\lam, \mu)).$ \begin{theorem}\label{thm:charlin} Let $Q(\lam, \mu) = \lam^2A_{20} + \mu^2A_{02}+\lam\mu A_{11}+\lam A_{10}+\mu A_{01}+A_{00}$ be a quadratic two-parameter matrix polynomial with real or complex coefficient matrices of order $n\times n$, and $v\in\C^3.$ Then a linear two-parameter matrix polynomial $L(\lam, \mu) \in \L(Q(\lam,\mu))$ corresponding to the ansatz vector $v$ is of the form $L(\lam, \mu)= \lam \widehat{A}_1 + \mu \widehat{A}_2 + \widehat{A}_3$ where \beano \widehat{A}_1 &=& \bmatrix{v\otimes A_{20} & -Y_1+v\otimes A_{11} &  -Z_1+v\otimes A_{10}}\\ \widehat{A}_2 &=& \bmatrix{Y_1 & v\otimes A_{02} & -Z_2+v\otimes A_{01} }\\ \widehat{A}_3 &=& \bmatrix{Z_1 & Z_2 & v\otimes A_{00}}\eeano where $Y_1=\bmatrix{Y_{11} \\ Y_{21} \\ Y_{31}}, Z_1=\bmatrix{Z_{11} \\ Z_{21} \\ Z_{31}}, Z_2=\bmatrix{Z_{12} \\ Z_{22} \\ Z_{32}}\in \C^{3n\times n}$ are arbitrary.

\end{theorem}

\noin\pf Let $\mathcal{M}: \L(Q(\lam, \mu)) \rightarrow \mathcal{V}_Q$ be a multiplicative map defined by $L(\lam, \mu)\mapsto L(\lam, \mu)(\Lambda\otimes I_n).$ Its easy to see that $\mathcal{M}$ is linear. First we show that $\mathcal{M}$ is surjective. Let $v\otimes Q(\lam, \mu)$ be an arbitrary element of $\mathcal{V}_Q.$ Construct $L(\lam, \mu)= \lam \widehat{A}_1 + \mu \widehat{A}_2 + \widehat{A}_3$ where \beano \widehat{A}_1 &=& \bmatrix{v\otimes A_{20} & v\otimes A_{11} &  v\otimes A_{10}}\\ \widehat{A}_2 &=& \bmatrix{0 & v\otimes A_{02} & v\otimes A_{01} }\\ \widehat{A}_3 &=& \bmatrix{0 & 0 & v\otimes A_{00}}.\eeano Then obviously we have $\widehat{A}_1 \boxplus \widehat{A}_2 \boxplus \widehat{A}_3 =v\otimes \bmatrix{A_{20} & A_{02} & A_{11} & A_{10} & A_{01} & A_{00}},$ so by Lemma \ref{lem:sec2} $L(\lam, \mu)$ is an $\mathcal{M}$-pre-image of $v\otimes Q(\lam, \mu).$ The set of all $\mathcal{M}$-preimages of $v\otimes Q(\lam, \mu)$ is $L(\lam, \mu)+ \mbox{Ker}\mathcal{M}$, so all that remains is to compute $\mbox{Ker}\mathcal{M}$. Further by Lemma \ref{lem:sec2} $\mbox{Ker}\mathcal{M}$ contains $L(\lam, \mu)= \lam \widehat{A}_1 + \mu \widehat{A}_2 + \widehat{A}_3$ that satisfies $\widehat{A}_1 \boxplus \widehat{A}_2 \boxplus \widehat{A}_3 = 0.$ The definition of ``box-addition" implies that  $\widehat{A}_1, \widehat{A}_2, \widehat{A}_3 $ are of the following form \beano \widehat{A}_1 &=& \bmatrix{0 & -Y_1 &  -Z_1}\\ \widehat{A}_2 &=& \bmatrix{Y_1 & 0 & 0 }\\ \widehat{A}_3 &=& \bmatrix{Z_1 & Z_2 & 0}\eeano where $Y_1, Z_1, Z_2 \in \C^{3n\times n}$ are arbitrary. This completes the proof.

\begin{example}\label{exp2} In Example \ref{exp1} we achieve the linear two-parameter polynomial $L(\lam,\mu)\in\L(Q(\lam,\mu))$ by choosing $$Y_1:=\bmatrix{-A_{20} \\ -A_{00}+A_{11} \\ A_{02}}, Z_1:=\bmatrix{-A_{01} \\ A_{10} \\ -I + 2A_{10}}, Z_2:=\bmatrix{0 \\ A_{01} \\ A_{01}}.$$ The \textit{standard} linearization $L(\lam,\mu)\in\L(Q(\lam,\mu))$ is achieved by choosing $$Y_1:=\bmatrix{0 \\ 0 \\ 0}, Z_1:=\bmatrix{A_{10} \\ 0 \\ -I}, Z_2:=\bmatrix{A_{01} \\ -I \\ 0}.$$ \end{example}
\begin{corollary} Dimension of $\L(Q(\lam,\mu))= 9n^2 + 3.$ \end{corollary}

\begin{remark} For quadratic one-parameter matrix polynomial \begin{equation}\label{def:pol} P(\lam) = \lam^2 A_{20} + \lam A_{10} + A_{00}, A_{i0}\in\C^{n\times n}, i=0,1,2, \end{equation} a vector space $\L_1(P)$ of matrix pencils of the form $L(\lam)=X+\lam Y \in\C^{2n\times 2n}$ is obtained in \cite{mackey3}. Setting $\mu=0$ in $Q(\lam,\mu)$ we have $Q(\lam,0)=P(\lam).$ Then from the constructions of linear two-parameter polynomials given in Theorem \ref{thm:charlin} it is easy to check that $\L(Q(\lam,\mu)) = \L_1(P).$ In fact, if $\mu =0$ then $\L(Q(\lam, \mu))$ contains matrix pencils $L(\lam)=\lam \widehat{A}_1 + \widehat{A}_3 \in \C^{2n\times 2n}$ where $$\widehat{A}_1 = \bmatrix{v\otimes A_{20}  &  -Z_1+v\otimes A_{10}}, \widehat{A}_3 = \bmatrix{Z_1 & v\otimes A_{00}},$$ $v\in\C^2$ and $Z_1\in\C^{2n\times n}$ is arbitrary. Thus we obtain the same vector space of matrix pencils obtained in \cite{mackey3} for a given quadratic one-parameter matrix polynomial $Q(\lam, 0)=\lam^2 A_{20} + \lam A_{10} + A_{00} = P(\lam).$
\end{remark}



\subsection{Construction of linearizations}

It is not very clear that whether all linear two-parameter matrix polynomials in the space $\L(Q(\lam,\mu))$ are linearizations of $Q(\lam, \mu).$ For example, consider any $L(\lam,\mu)\in \L(Q(\lam,\mu))$ corresponding to \textit{ansatz} vector $v = 0$. Thus given a quadratic two-parameter matrix polynomial $Q(\lam,\mu)$ we need to identify which $L(\lam,\mu)$ in $\L(Q(\lam,\mu))$ are linearizations.


 We begin with a result concerning the special case of the \textit{ansatz} vector $v=\alpha e_1$ where $e_1=\bmatrix{1&0&0}^T$ and $0\neq\alpha\in\C.$



\begin{theorem}\label{thm:linearization} Let $Q(\lam, \mu)=\lam^2 A_{20} + \lam\mu A_{11} + \mu^2 A_{02} + \lam A_{10} + \mu A_{01} + A_{00}$ be a quadratic two-parameter matrix polynomial with real or complex coefficient matrices of order $n\times n$. Suppose  $L(\lam, \mu)=\lam\widehat{A}_1 + \mu\widehat{A}_2 + \widehat{A}_3 \in \L(Q(\lam, \mu))$ with respect to the ansatz vector $v=\alpha e_1\in\C^3,$  where \beano \widehat{A}_1 &=& \bmatrix{\alpha e_1\otimes A_{20} & -Y_1+\alpha e_1\otimes A_{11} &  -Z_1+\alpha e_1\otimes A_{10}}\\ \widehat{A}_2 &=& \bmatrix{Y_1 & \alpha e_1\otimes A_{02} & -Z_2+\alpha e_1\otimes A_{01} }\\ \widehat{A}_3 &=& \bmatrix{Z_1 & Z_2 & \alpha e_1\otimes A_{00}},\eeano $Y_1=\bmatrix{Y_{11} \\ 0 \\ 0}, Z_1=\bmatrix{Z_{11} \\ Z_{21} \\ Z_{31}}, Z_2=\bmatrix{Z_{12} \\ Z_{22} \\ Z_{32}}\in \C^{3n\times n}, \,\, \mbox{det}\bmatrix{Z_{21} & Z_{22} \\ Z_{31} & Z_{32}} \neq 0.$ Then $L(\lam,\mu)$ is a linearization of $Q(\lam, \mu).$
\end{theorem}

\noin\pf By Theorem \ref{thm:charlin}, any linear two-parameter matrix polynomial $L(\lam, \mu) = \lam\widehat{A}_1 + \mu\widehat{A}_2 + \widehat{A}_3 \in \L(Q(\lam,\mu))$ corresponding to the \textit{ansatz} vector $v=\alpha e_1$ is of the form \begin{small} \beano L(\lam,\mu) = && \lam \bmatrix{\alpha A_{20} & -Y_{11}+\alpha A_{11} & -Z_{11}+\alpha A_{10} \\ 0 & -Y_{21} & -Z_{21} \\ 0 & -Y_{31} & -Z_{31}} + \mu \bmatrix{Y_{11} & \alpha A_{02} & -Z_{12}+\alpha A_{01} \\ Y_{21} & 0 & -Z_{22} \\ Y_{31} & 0 & -Z_{32}} \\ && + \bmatrix{Z_{11} & Z_{12} & \alpha A_{00} \\ Z_{21} & Z_{22} & 0 \\ Z_{31} & Z_{32} & 0}. \eeano \end{small}
Thus we have $$L(\lam,\mu) = \bmatrix{W_1(\lam, \mu) & W_2(\lam, \mu) & W_3(\lam, \mu) \\ \mu Y_{21}+Z_{21} & -\lam Y_{21}+Z_{22} & -\lam Z_{21}-\mu Z_{22} \\ \mu Y_{31}+Z_{31} & -\lam Y_{31}+Z_{32} & -\lam Z_{31}-\mu Z_{32}}$$ where $W_1(\lam, \mu) = \alpha \lam A_{20}+\mu Y_{11}+Z_{11}, W_2(\lam, \mu) = \alpha \mu A_{02}+\lam \alpha A_{11}-\lam Y_{11}+Z_{12}$ and $W_3(\lam, \mu) = \alpha\lam A_{10}-\lam Z_{11}+\alpha \mu A_{01}-\mu Z_{12}+\alpha A_{00}.$

 Define $$E(\lam,\mu) = \bmatrix{\frac{\lam}{\alpha} I & I & 0 \\ \frac{\mu}{\alpha} I & 0 & I \\ \frac{1}{\alpha} I & 0 &0 }.$$ Consequently, we have $$L(\lam,\mu)E(\lam,\mu) = \bmatrix{Q(\lam,\mu) & W_1(\lam,\mu) & W_2(\lam,\mu) \\ 0 & \mu Y_{21}+Z_{21} & -\lam Y_{21}+Z_{22} \\ 0 & \mu Y_{31}+Z_{31} & -\lam Y_{31}+Z_{32}}.$$ Setting $Y_{21}=0=Y_{31}$ we have $L(\lam,\mu)E(\lam,\mu) = \bmatrix{Q(\lam,\mu) & W(\lam,\mu) \\ 0 & Z}$ where $W(\lam,\mu)=\bmatrix{W_1(\lam,\mu) & W_2(\lam,\mu)}\in\C^{n\times 2n}, Z=\bmatrix{Z_{21} & Z_{22} \\ Z_{31} & Z_{32}}\in\C^{2n\times 2n}.$

Since $Z$ is nonsingular, we define $$F(\lam,\mu) = \bmatrix{I & -W(\lam,\mu)Z^{-1} \\ 0 & Z^{-1}}.$$ Then we have $$F(\lam,\mu)L(\lam,\mu)E(\lam,\mu)= \bmatrix{ Q(\lam,\mu) & 0 \\ 0 & I_{2n}}.$$ Note that both $E(\lam, \mu)$ and $F(\lam, \mu)$ are unimodular polynomials. Hence we have $\mbox{det}L(\lam,\mu)=\gamma\mbox{det}Q(\lam,\mu)$ for some nonzero $\gamma\in\C.$ Thus $L(\lam,\mu)$ is a linearization. This completes the proof.


 Let $Q(\lam,\mu)$ quadratic two-parameter matrix polynomial and $L(\lam, \mu)\in\L(Q(\lam, \mu))$ corresponding to an \textit{ansatz} vector $0\neq v\in\C^3.$ Then the following is a procedure for determining a set of linearizations of $Q(\lam, \mu).$\\

\noin\textbf{Procedure to determine linearizations in $\L(Q(\lam,\mu))$:}

\begin{enumerate}
\item Suppose $Q(\lam,\mu)$ is a quadratic two-parameter matrix polynomial and $L(\lam,\mu) = \lam¸\widehat{A}_1 +\mu \widehat{A}_2 + \widehat{A}_3 \in \L(Q(\lam,\mu))$ corresponding to \textit{ansatz} vector $v\in\C^3$ i.e. $L(\lam,\mu) (\Lambda\otimes I_n) = v\otimes Q(\lam,\mu).$

     \item Select any nonsingular matrix $M=\bmatrix{m_{11}&m_{12}&m_{13}\\m_{21}&m_{22}&m_{23}\\m_{31}&m_{32}&m_{33}}$ such that $Mv = \alpha e_1\in \C^3, \alpha\neq 0.$ A list of nonsingular matrices $M$ depending on the entries of $v$ is given in the Appendix.

     \item Apply the corresponding block-transformation $M \otimes I_n$ to $L(\lam,\mu).$ Then we have $\widetilde{L}(\lam, \mu)= (M\otimes I_n)L(\lam,\mu)= \lam \widetilde{A_1} + \mu \widetilde{A_2} + \widetilde{A_3}$ such that \beano \widetilde{A_1} &=& \bmatrix{\alpha e_1\otimes A_{20} & -\widetilde{Y_1}+\alpha e_1\otimes A_{11} &  -\widetilde{Z_1}+\alpha e_1\otimes A_{10}}\\ \widetilde{A_2} &=& \bmatrix{\widetilde{Y_1} & \alpha e_1\otimes A_{02} & -\widetilde{Z_2}+\alpha e_1\otimes A_{01} }\\ \widetilde{A_3} &=& \bmatrix{\widetilde{Z_1} & \widetilde{Z_2} & \alpha e_1\otimes A_{00}}\eeano where \beano \widetilde{Y_1} &=& (M\otimes I_n)Y_1=(M\otimes I_n)\bmatrix{Y_{11}\\ 0\\ 0} = \bmatrix{m_{11}Y_{11} \\ m_{21}Y_{11} \\ m_{31}Y_{11}} \\ \widetilde{Z_1} &=& (M\otimes I_n)Z_1=(M\otimes I_n)\bmatrix{Z_{11}\\ Z_{21}\\ Z_{31}} =\bmatrix{m_{11}Z_{11} +m_{12}Z_{21} + m_{13}Z_{31} \\ m_{21}Z_{11} +m_{22}Z_{21} + m_{23}Z_{31} \\ m_{31}Z_{11} +m_{32}Z_{21} + m_{33}Z_{31}} \\  \widetilde{Z_2}&=&(M\otimes I_n)Z_2 = (M\otimes I_n)\bmatrix{Z_{12}\\ Z_{22}\\ Z_{32}}=\bmatrix{m_{11}Z_{12} +m_{12}Z_{22} + m_{13}Z_{32} \\ m_{21}Z_{12} +m_{22}Z_{22} + m_{23}Z_{32} \\ m_{31}Z_{12} +m_{32}Z_{22} + m_{33}Z_{32}}\eeano are arbitrary.

         \item For $\widetilde{L}(\lam, \mu)$ to be linearization we need to choose $Y_1, Z_1, Z_2$ as follows. If $m_{21}=m_{31}=0$ then choose $Y_{11}$ arbitrary; otherwise choose $Y_{11}=0.$  Further we need to choose $Z_1=\bmatrix{Z_{11}\\Z_{21}\\Z_{31}}, Z_2=\bmatrix{Z_{12}\\Z_{22}\\Z_{32}}$ for which \begin{equation}\label{eqn:cond}\mbox{det}\bmatrix{m_{21}Z_{11} +m_{22}Z_{21} + m_{23}Z_{31} & m_{21}Z_{12} +m_{22}Z_{22} + m_{23}Z_{32}  \\ m_{31}Z_{11} +m_{32}Z_{21} + m_{33}Z_{31} & m_{31}Z_{12} +m_{32}Z_{22} + m_{33}Z_{32}}\neq 0.\end{equation} From the construction of $M$ given in the Appendix it is easy to check that we can always choose suitable $Z_1, Z_2$ for which the condition (\ref{eqn:cond}) is satisfied.



\end{enumerate}


\section{Linearization of two-parameter quadratic eigenvalue problem}


%
The quadratic two-parameter eigenvalue problem is concerned with finding a pair $(\lam,\mu)\in\C\times\C$ and nonzero vectors $x_i\in\C^{n_i}$ for which
\be\label{defn:qevp2} Q_i(\lam,\mu)x_i = 0, i=1,2, \ee
where \be\label{def:q2pp} Q_i(\lam,\mu)=A_i + \lam B_i + \mu C_i + \lam^2 D_i + \lam\mu E_i + \mu^2 F_i, \ee $A_i, B_i, \hdots, F_i\in\C^{n_i\times n_i}.$
The pair $(\lam,\mu)$ is called an eigenvalue of (\ref{defn:qevp2}) and $x_1\otimes x_2$ is called the corresponding eigenvector. The spectrum of a quadratic two-parameter eigenvalue problem is the set \be\sigma_Q :=\left\{ (\lam,\mu)\in\C\times\C : \mbox{det}Q_i(\lam,\mu)=0, i=1,2\right\}.\ee In the generic case, we observe that (\ref{defn:qevp2}) has $4n_1n_2$ eigenvalues by using the following theorem.

\begin{theorem}\label{thm:Bezout}(Bezout's Theorem, \cite{CoxLiOs97}) Let $f(x,y)=g(x,y)=0$ be a system of two polynomial equations in two unknowns. If it has only finitely many common complex zeros $(x,y)\in\C\times\C,$ then the number of those zeros is at most $\mbox{degree}(f)\cdot\mbox{degree}(g).$ \end{theorem}


The usual approach to solving (\ref{defn:qevp2}) is to linearize it as a two-parameter eigenvalue problem given by \begin{equation}\label{defn:qlin2} \left\{
                     \begin{array}{ll}
                       L_1(\lam,\mu)w_1 = (A^{(1)}) + \lam B^{(1)} + \mu C^{(1)})w_1=0 & \hbox{} \\
                       L_2(\lam,\mu)w_2 = (A^{(2)}) + \lam B^{(2)} + \mu C^{(2)})w_2=0 & \hbox{}
                     \end{array}
                   \right.
   \end{equation} where $A^{(i)}, B^{(i)}, C^{(i)} \in \C^{m_i\times m_i}, m_i\geq 2n_i, i=1,2,$ and $w_i = \Lambda \otimes x_i. $ A pair $(\lam,\mu)$ is called an eigenvalue of (\ref{defn:qlin2}) if $L_i(\lam, \mu)w_i = 0$ for a nonzero vector $w_i$ for $i =1, 2$, and $w_1\otimes w_2$ is the corresponding eigenvector. Thus the spectrum of the linearized two-parameter eigenvalue problem is given by \be\sigma_L :=\left\{ (\lam,\mu)\in\C\times\C : \mbox{det}L_i(\lam,\mu)=0, i=1,2\right\}.\ee Therefore, in the generic case, the problem (\ref{defn:qlin2}) has $m_1m_2\geq 4n_1n_2$ eigenvalues.

A standard approach to solve a two-parameter eigenvalue problem (\ref{defn:qlin2})
is by converting it into a coupled generalized eigenvalue problem given by \be \Delta_1 z=\lam \Delta_0 z \,\, \mbox{and} \,\, \Delta_2 z=\mu \Delta_0 z\ee where $z=w_1\otimes w_2$ and \beano \Delta_0 &=& B^{(1)} \otimes C^{(2)} - C^{(1)}\otimes B^{(2)} \\ \Delta_1 &=& C^{(1)}\otimes A^{(2)} - A^{(1)}\otimes C^{(2)} \\ \Delta_2 &=& A^{(1)}\otimes B^{(2)} - B^{(1)}\otimes A^{(2)}.\eeano
The two-parameter eigenvalue problem is called singular (resp. nonsingular) if $\Delta_0$ is singular (resp. nonsingular), see \cite{mupl08}.

As mentioned earlier, we are interested in finding linear two-parameter polynomials $L_i(\lam,\mu)$ for a given quadratic two-parameter eigenvalue problem (\ref{defn:qevp2}) such that $\sigma_Q=\sigma_L.$ Thus we have the following definition.

\begin{definition} Let (\ref{defn:qevp2})
be a quadratic two-parameter eigenvalue problem. A two-parameter eigenvalue problem (\ref{defn:qlin2})
is said to be a linearization of (\ref{defn:qevp2}) if $L_i(\lam, \mu)$ is a linearization of $Q_i(\lam, \mu).$ \end{definition}

 Thus if we consider a linearization of a quadratic two-parameter eigenvalue problem then $\sigma_Q = \sigma_L$ is guaranteed. It is also easy to observe that $x_1\otimes x_2$ is an eigenvector corresponding to an eigenvalue $(\lam, \mu)$ of a quadratic two-parameter eigenvalue problem if and only if $w_1\otimes w_2 $ is an eigenvector corresponding to the eigenvalue $(\lam, \mu)$ of the linearization.

%

  Making use of the construction of linearizations for a two-parameter quadratic matrix polynomial described in section 2, we construct linearizations for a quadratic two-parameter eigenvalue problem.

\begin{theorem}\label{Thm:q2pp} Let (\ref{defn:qevp2})
be a quadratic two-parameter eigenvalue problem. A class of linearizations of (\ref{defn:qevp2}) is given by $$L_i(\lam,\mu)w_i= (A^{(i)} + \lam B^{(i)} + \mu C^{(i)})w_i=0, w_i=\Lambda\otimes x_i, i=1,2,$$ where \beano A^{(i)} &=& \bmatrix{Z_1^{(i)} & Z_2^{(i)} & \alpha_ie_1\otimes A_i }, \\ B^{(i)} &=& \bmatrix{\alpha_ie_1\otimes D_i & -Y_1^{(i)} + \alpha_ie_1\otimes E_i & -Z_1^{(i)}+\alpha_ie_1\otimes B_i}, \\ C^{(i)} &=& \bmatrix{Y_1^{(i)} & \alpha_ie_1\otimes F_i & -Z_2^{(i)}+\alpha_ie_1\otimes C_i},\eeano $\alpha_i\neq 0, Y_1^{(i)}=\bmatrix{Y_{11}^{(i)} \\ 0 \\ 0}, Z_1^{(i)} = \bmatrix{Z_{11}^{(i)} \\Z_{21}^{(i)} \\ Z_{31}^{(i)}  }, Z_2^{(i)} = \bmatrix{Z_{12}^{(i)} \\Z_{22}^{(i)} \\ Z_{32}^{(i)}  }  \in \K^{3n\times n}, \mbox{det}\bmatrix{Z_{21}^{(i)} & Z_{22}^{(i)} \\ Z_{31}^{(i)} & Z_{32}^{(i)}} \neq 0.$
\end{theorem}

\pf Consider the linearizations $L_i(\lam,\mu)= A^{(i)} + \lam B^{(i)} + \mu C^{(i)}$ of $Q_i(\lam,\mu)$ associated with \textit{ansatz} vector $0\neq \alpha_ie_1\in \C^3, i=1,2$ given by Theorem \ref{thm:linearization}. This completes the proof.


Now we show that the linearizations for a quadratic two-parameter eigenvalue problem described in Theorem \ref{Thm:q2pp} are singular linearizations. The following theorem plays an important role in the sequel.

\begin{theorem}\label{Thm:det} The determinant of a block-triangular matrix is the product of the determinants of the diagonal blocks.\end{theorem}

Now we have the following result.

\begin{theorem}
The linearizations for (\ref{defn:qevp2}) derived in Theorem \ref{Thm:q2pp} are singular linearizations.

\end{theorem}


\pf Consider the linearizations $L_i(\lam,\mu)w_i=  (A^{(i)} + \lam B^{(i)} + \mu C^{(i)})w_i = 0, i=1,2$ of $Q_i(\lam,\mu)$ where \beano B^{(i)} = \bmatrix{\alpha_iD_i & -Y_{11}^{(i)} + \alpha_i E_i & -Z_{11}^{(i)} + \alpha_i B_i \\ 0 & 0 & -Z_{21}^{(i)} \\ 0 & 0 & -Z_{31}^{(i)} }, C^{(i)} = \bmatrix{Y_{11}^{(i)}& \alpha_i F_i & -Z_{12}^{(i)} + \alpha_i C_i \\ 0 & 0 & -Z_{22}^{(i)} \\ 0 & 0 & -Z_{32}^{(i)} }. \eeano Consequently we have \beano \Delta_0 &=& B^{(1)}\otimes C^{(2)} - C^{(1)}\otimes B^{(2)} \\ &=& \bmatrix{\alpha_1D_1\otimes C^{(2)} & (-Y_{11}^{(1)} + \alpha_1 E_1)\otimes C^{(2)}  & (-Z_{11}^{(1)} + \alpha_1 B_1)\otimes C^{(2)}  \\ 0 & 0 & -Z_{21}^{(1)}\otimes C^{(2)}  \\ 0 & 0 & -Z_{31}^{(1)}\otimes C^{(2)}  } \\ && - \bmatrix{Y_{11}^{(1)}\otimes B^{(2)} & \alpha_1 F_1\otimes B^{(2)} & (-Z_{12}^{(1)} + \alpha_1 C_1)\otimes B^{(2)} \\ 0 & 0 & -Z_{22}^{(1)}\otimes B^{(2)}\\ 0 & 0 & -Z_{32}^{(1)}\otimes B^{(2)} }. \eeano Observe that $\Delta_0$ is a block-triangular matrix with one of the diagonal blocks is $0.$ Hence by Theorem \ref{Thm:det} we have $\mbox{det}\Delta_0 =0.$ This completes the proof.



\begin{remark} Note that given a quadratic two-parameter eigenvalue problem (\ref{defn:qevp2})
we choose linearizations $L_i(\lam,\mu)$ of $Q_i(\lam,\mu)$ associated with the \textit{ansatz} vector $0\neq \alpha_ie_1\in\C^3$, and constructed linearizations $L_i(\lam,\mu)w_i=0, w_i=\Lambda\otimes x_i$ of (\ref{defn:qevp2}). However, we can derive a large class of singular linearizations by choosing linearizations $L_i(\lam,\mu)$ of $Q_i(\lam,\mu)$ associated with \textit{ansatz} vector $0\neq v_i\in\C^3$ described in section $2$. \end{remark}

\section{Conclusions} Given a quadratic two-parameter matrix polynomial $Q(\lam,\mu)$, we construct a vector space of linear two-parameter matrix polynomials and identify a set of linearizations of $Q(\lam,\mu).$ We also describe construction of each of these linearizations. Finally, using these linearizations we determine a class of singular linearizations for a quadratic two-parameter eigenvalue problem.\\

\noin\textbf{Acknowledgement} The author wishes to thank Bor Plestenjak and Rafikul Alam for their stimulating comments that have significantly improved the quality of the results.  Thanks are due to the referees for their constructive remarks.\\\\

\noin \textbf{Appendix}

Let $0\neq \alpha\in\C$ and $e_1=\bmatrix{1 \\0 \\0}.$ Given a vector $v=\bmatrix{a\\b\\c}\in\C^3$ we can always pick a nonsingular matrix $M\in\C^{3\times 3}$ for which $Mv=\alpha e_1$ as follows.
$$M=\left\{
    \begin{array}{ll}
      \small{\bmatrix{\alpha/a & 0 & 0 \\ 1/a & -1/b & 0 \\ 1/a & 0 & -1/c}}, & \hbox{if $a\neq 0, b\neq 0, c\neq 0$} \\
      \small{\bmatrix{0 & \alpha/b & 0 \\ 0 & -1/b & 1/c \\ 1 & 0 & 0}}, & \hbox{if $a=0, b\neq 0, c\neq 0$} \\
      \small{\bmatrix{1 & 1 & \alpha/c \\ 1 & 1 & 0 \\ 0 & 1 & 0}}, & \hbox{if $a= 0, b= 0, c\neq 0$}\\
\small{\bmatrix{\alpha/a & 0 & 0 \\ 0 & 1 & 0 \\ -1/a & 0 & 1/c}}, & \hbox{if $a\neq 0, b= 0, c\neq 0$}\\
\small{\bmatrix{\alpha/a & 0 & 0 \\ 0 & 1 & 0 \\ 0 & 1 & 1}}, & \hbox{if $a\neq 0, b= 0, c= 0$}\\
\small{\bmatrix{\alpha/a & 0 & 1 \\ 1/a & -1/b & 1 \\ -1/a & 1/b & 0}}, & \hbox{if $a\neq 0, b\neq 0, c= 0$}\\
\small{\bmatrix{1 & \alpha/b & 0 \\ 1 & 0 & 0 \\ 1 & 0 & 1}}, & \hbox{if $a= 0, b\neq 0, c= 0$}\\
\small{\bmatrix{\alpha/a & 0 & 0 \\ 1/a & 0 & -1/c \\ 0 & 1 & 0}}, & \hbox{if $a\neq 0, b= 0, c\neq 0.$}
    \end{array}
  \right.$$

%
%
%
%
%
%
%

\end{document}